\documentclass[11pt]{article}

\usepackage{amsmath}
\usepackage{amsthm}

\usepackage[T1]{fontenc}

\usepackage{color}
\usepackage{graphicx}
\usepackage{amssymb}
\usepackage{lscape}
\usepackage{enumerate}
\usepackage{mathrsfs}
\usepackage{booktabs}
\usepackage{wrapfig}
\usepackage{float}
\usepackage{floatflt}
\usepackage[caption = false]{subfig}

\def\D{{\mathcal{ D}}}
\def\V{\mathcal{V}}
\def\M{\mathcal{M}}

\def\1{{\sf {1}}}
\def\0{{\sf {0}}}

\def\n{{\sf {n}}}

\def\demo{{\noindent\bf Proof.}\hskip 0.3truecm}
\def\BBox{\kern  -0.2cm\hbox{\vrule width 0.15cm height 0.3cm}}

\oddsidemargin .cm \evensidemargin 1.cm \textwidth=16.5cm
\include {mak}
\begin{document}
\newtheorem{propo}{Proposition}[section]
\newtheorem{lemma}[propo]{Lemma}
\newtheorem{observation}[propo]{Observation}
\newtheorem{theorem}[propo]{Theorem}
\newtheorem{corollary}[propo]{Corollary}
\newtheorem{definition}[propo]{Definition}
\newtheorem{example}[propo]{Example}

\def\RR{\mathbb{R}}
\def\NN{\mathbb{N}}
\def\1{{\sf 1}}
\def\n{{\sf n}}
\def\so{{\sf supp}}
\def\proof{{\noindent  Proof.}\hskip 0.3truecm}
\def\BBox{\kern  -0.2cm\hbox{\vrule width 0.15cm height 0.3cm}}
\def\endproof{\hspace{.20in} \BBox\vspace{.20in}}
\def\L{{\mathscr  L}}\def\P{{\mathscr  P}}
\def\C{{\mathcal C}}
\def\K{{\mathscr K}}
\def\E{{\mathcal E}}
\def\F{{\mathcal F}}
\def\G{\mathcal{G}}
\def\P{{\mathscr P}}
\def\L{{\mathcal L}}
\def\C{{\mathcal C}}
\def\D{{\mathcal D}}
\def\V{{\mathcal V}}
\def\F{{\mathcal F}}
\def\E{{\mathcal E}}
\def\M{{\mathcal M}}
\def\O{{\mathcal O}}
\def\N{{\mathcal N}}
\def\R{{\mathcal R}}
\def\mG{{\widetilde \G}}
\def\mP{{\widetilde \P}}
\def\mR{{\widetilde \R}}
\def\mkG{{\widetilde G}}
\def\mkP{{\widetilde P}}
\def\mkR{{\widetilde R}}
\def\dr{{N}}
\def\mPP{{\widetilde {\sf P}}}
\def\mRR{{\widetilde {\sf R}}}
\oddsidemargin 16.5mm
\evensidemargin 16.5mm

\title{Effective resistances and Kirchhoff Index in subdivision networks}
\author{{\small \'Angeles Carmona, Margarida Mitjana, Enric Mons\'o  }\\
\\
{\small \em Departament de Matem\`atiques}\\
{\small Universitat Polit\`ecnica de Catalunya. Spain}\\
}
\date{}
\maketitle

\begin{abstract}
 We define a  {\it subdivision network} $\Gamma^S$ of a given network $\Gamma,$  by inserting a new vertex in every edge, so that each edge  is
replaced by two new edges with  conductances that fulfill electrical conditions on the new network.  In this work, we firstly obtain an expression for the Green kernel of the subdivision network in terms of the Green kernel of the base network. Moreover, we also obtain the effective resistance and the Kirchhoff index of the subdivision network in terms of the corresponding parameters on the base network.   Finally, as an example, we carry out the computations in the case of a wheel.
\end{abstract}


\section{Introduction}
Many recent papers are devoted to the study of different parameters of the subdivison graphs. For instance, Chen in  \cite{Che10}, obtained a formula for the effective resistances of the subdivision graph in terms of the effective resistances of the original graph by using some nice sum rules. The Kirchhoff index of the subdivision graph is considered in different works under several hypothesis such as regular graphs in \cite{GaLuLi12}, bipartite graphs in \cite{SuWaZhBu15}, or operations between graphs that involve the subdivision concept as well, see \cite{BuYaZhZh14} for instance.

In \cite{Ya14}, the author extends the previous results to general graphs and computes the Kirchhoff index of subdivision graph  in terms of the Kirchhoff index, the multiplicative degree--Kirchhoff index,
the additive degree--Kirchhoff index, the number of vertices, and the number of edges of $\Gamma$. Simultaneously, Sun {\it et alt.}   gave the formulae for the Kirchhoff index in terms of a $\{1\}$--inverse of the combinatorial Laplacian, see \cite{SuWaZhBu15}.

In the present paper, we  introduce the   subdivision of a network. Our approach consists in interpreting a   network as an electric circuit, and hence each edge has got assigned a positive number that corresponds with the conductance of a wire connecting two nodes,  its inverse is the resistance. When we perform the subdivision operation we interpret that we introduce a rheostat in every edge, that is a device  that may change the resistance without opening the circuit in which it is connected. Thus, we decompose each edge in two new edges taking into account electrical compatibility of the circuit, specifically, the series sum rule for resistances. As a consequence, we would get that after the subdivision process, the effective resistance between any pair of old vertices  should remain unchanged.

Our methodology comes from discrete Potential theory and hence, we express all the parameters in terms of the Green kernel of the network. That can be seen in matrix terms, as the computation of the Group inverse of the combinatorial Laplacian of a subdivision network in terms of the Group inverse of the combinatorial Laplacian of the base network. 

In Section \ref{Sec.PoissonProblem} we first  obtain a solution of the Poisson problem in the subdivision network in terms of the solution of an appropriate Poisson problem on the base network and hence we compute the Green kernel of the subdivision network. Next, we give an expression for the  effective resistance between any pair of vertices of the subdivision network and its corresponding Kirchhoff index. For all the results we compare ours with the previously known for the case of graphs, and in particular for $k$--regular graphs. 

 The last section contains  the expressions for the Green kernel, the effective resistance and the Kirchhoff index of the subdivision network of a wheel as an illustration of the obtained results.

We end the present section by introducing the basic notation  and results.

In the whole work, a {\it network} is the triplet
$\Gamma=(V,E,c)$ where  $(V,E)$ stands for a finite and connected graph, without loops
nor multiple edges; and $c\colon V\times
V\longrightarrow [0,+\infty)$ is a symmetric function called {\it
conductance} satisfying $c(x,y)>0$ iff $x\sim y$ which means that $\{x,y\}\in E.$ Let $n$ be the number of nodes and  $m$ the
number of edges.

 On the other hand, $\mathcal{C}(V)$ is the set of real functions
on $V$. For any vertex $x\in V,$  $\varepsilon_x \in  \mathcal{C}(V)$ is the Dirac function at $x$ and  $k\in \mathcal{C}(V)$ defined as
$k(x)=\sum\limits_{y\in V}c(x,y),$  is the {\it degree} of $x$. The
standard inner product in $\mathcal{C}(V)$ is denoted by
$\left\langle\cdot,\cdot\right\rangle$; that is, if $u,v \in \mathcal{C}(V)$ then, $\left\langle u,v\right\rangle=\sum\limits_{x\in V}u(x)v(x).$

The {\it Laplacian} of $\Gamma$ is the linear operator
$\mathcal{L}\colon \mathcal{C}(V)\longrightarrow \mathcal{C}(V)$
defined, for each $u\in \mathcal{C}(V)$ and $x\in V$ as
$$\mathcal{L}(u)(x)=\sum\limits_{y \in V}c(x,y)\big(u(x)-u(y)\big).$$
 The Laplacian of $ \Gamma$ is a self--adjoint and positive
semi--definite operator. Moreover, $\L(u)=0$ iff $u$ is constant and hence,  ${\sf
ker}({\mathcal L})={\rm span}(\1)$, where $\1\in\C(V)$ is the function that assigns $1$ to any vertex. Therefore,  ${\mathcal L}$  defines  an
isomorphism on  ${\sf ker}({\mathcal L})^\bot$.

A  {\it Poisson problem} consists in, given $f\in\C(V)$,  finding $u\in\C(V)$ such that
\begin{equation}
\label{poisson}
\L(u)=f \hspace{.25cm} \mbox{ on }\hspace{.25cm} V.
\end{equation}
From the above properties, the Poisson equation has solution iff $\langle f,\1\rangle=0$, and there is a unique solution of \eqref{poisson} satisfying $\langle u,\1\rangle=0.$

The operator that assigns to every function $f\in \mathcal{C}(V)$ the unique solution of $\mathcal{L}(u)=f-\displaystyle\langle f,\1\rangle\1,$ such that $\left\langle u,\1\right\rangle=0$,  is called the {\it Green operator} and it is denoted by $\G.$  The operator $\G$ is self--adjoint and positive semi--definite. Then the symmetric function $G,$
defined as $G(x,y)=\mathcal{G}(\varepsilon_y)(x)$ for every pair $x,y\in V,$ is called the {\it Green kernel} on $V.$  Moreover, $\displaystyle \mathcal{G}(f)(x)=\sum\limits_{x,y\in V} G(x,y)\,f(y)$ and $\langle\G(f),f\rangle=0$ iff $f=a\1$, $a\in \mathbb{R}$.

The relation between an integral operator and its associated kernel enables  us to characterize the Green  kernel for $\Gamma$ as solutions of suitable boundary value problems. For all $y\in V$, the function $G_y=G(\cdot,y)$ is cha\-rac\-terized by  equations
\begin{equation}\label{Greencharac}
{\mathcal L}(G_y)=\varepsilon_y-\dfrac{1}{n}\1 \hspace{.5cm}\mbox{ and }\hspace{.5cm}  \displaystyle\langle G_y,\1\rangle=0.\end{equation}
See  \cite{CaEnMi14-1} and references therein for more details. Notice that, if we label the vertices of $\Gamma$, both the Laplacian and the Green operator  can be interpreted as matrices and hence, the Green kernel can be identified as the Group Inverse of the combinatorial Laplacian.

The {\it effective resistance} between vertices $x$ and $y$ is defined as $R(x,y)=u(x)-u(y)$, where  $u\in \mathcal{C}(V)$ is any solution of the Poisson problem $\mathcal{L}(u)=\varepsilon_x-\varepsilon_y$. The effective resistance  can be interpreted as  the  voltage measured in  the nodes $x$ and $y$ when  a unitary current is applied between them.
Actually, $R$ defines a distance on $\Gamma$, usually referred as {\it resistive distance,} see \cite{BeCaEnGe08a,KlRa93}, and  gives a measure of how much two different nodes of a network are connected. Thus, the more well connected are vertices  $x$ and $y,$ the less is $R(x,y).$ Moreover, for any $x,y\in V$ the following relation holds
\begin{equation}\label{resistancegreen}
R(x,y)=G(x,x)+G(y,y)-
2G(x,y).
\end{equation}
The {\it Kirchhoff Index} ${\sf k}$ of a network $\Gamma$, also called its {\it total resistance},  is defined as
\begin{equation}
\label{D:kirchhoff} {\sf k}=\frac{1}{2}\sum\limits_{x,y\in V}R(x,y)=n \sum\limits_{x\in V}G(x,x)
\end{equation}
and gives a measure of the global connectivity of the network.  The Kirchhoff index is a descriptor of the structure of the network and exhibits many  interesting interpretations, see \cite{GhBoSa08, XiGu03}.


%
%
%


\section{The Poisson Problem on a Subdivision Network}
\label{Sec.PoissonProblem}
A {\it subdivision network} $\Gamma^S=(V^S, E^S, c^S)$ of a given network $\Gamma=(V, E,c),$ is
obtained by inserting a new vertex in every edge, so that each edge $\{x,y\}\in E$ is
replaced by two new edges, say $\{x, v_{xy}\}$  and $\{y,v_{xy}\}$ where $v_{xy}$ is the new inserted vertex. We denote by $V'$ the new vertex set assuming that, $v_{xy}=v_{yx}.$ Thus, $V^S=V\cup V',$ the order of the subdivision network is $n+m,$  whereas the size is $2m.$
Moreover, according to the well--known rule that express  the equivalent resistance of two resistors connected in series,  we define the conductance function $c^S\colon V^S\times V^S\longrightarrow [0,+\infty)$  by choosing, for every pair of adjacent vertices, non--null values $c^S(x,v_{xy})$ and $c^S(y,v_{xy})$ such that
\begin{equation}\label{conductances}
\frac{1}{c(x,y)}=\frac{1}{c^S(x,v_{xy})}+\frac{1}{c^S(y,v_{xy})}.\end{equation}
The definition of $c^S$ cannot be misunderstood as all the edges in
$E^S$ have both kind of vertices, one in $V$ and the other in $V'$. Hence, by the sake of simplicity, it will be denoted as $c.$  Moreover for each edge, there exist infinitely many different  choices of conductances fulfilling \eqref{conductances}, so that different choices will lead to different  subdivision networks.

Up to our knowledge, the only case that has been  studied in the literature, (\cite{Che10,GaLuLi12,SuWaZhBu15,Ya14}), is $c(x,y)=c(x,v_{xy})=c(y,v_{xy})=1$, that not fulfills the electrical compatibility  condition \eqref{conductances}. In the present work, and in order to compare with the known results, we will consider as a particular case $c(x,y)=1$ and $c(x,v_{xy})=c(y,v_{xy})=2$ and we call it {\it standard subdivision graph}.

Observe that $\Gamma^S$ is also a connected, finite, with no loops, nor multiple edges network.

If $\L^S$ denotes the combinatorial Laplacian of $\Gamma^S,$ then for any $u\in\C(V^S)$ we have that
$$
\begin{array}{rll}
\mathcal{L}^S(u)(x)=&\hspace{-.25cm}\displaystyle\sum\limits_{y\in V } c(x,v_{xy})\left(u(x)-u(v_{xy})\right),&\hbox{for any }x\in V;\\[2ex]
\mathcal{L}^S(u)(v_{xy})=&\hspace{-.25cm}c(x,v_{xy})\left(u(v_{xy})-u(x)\right)+c(y,v_{xy})\left(u(v_{xy})-u(y)\right),&\hbox{for any }v_{xy}\in V'.
\end{array}$$
%

The aim of this section is to obtain a solution of the Poisson problem in $\Gamma^S$ in terms of the solution of an appropriate Poisson problem
on $\Gamma.$

It is helpful for the sequel to define, for each pair $x,y\in V$ with $x\sim y$, the coefficient
$$\alpha(x,y)=\frac{c(x,v_{xy})}{c(x,v_{xy})+c(y,v_{xy})}=\frac{c(x,v_{xy})}{k(v_{xy})},$$
where $k(v_{xy})=c(x,v_{xy})+c(y,v_{xy})$, is the degree of $v_{xy}$ in $\Gamma^S$.
Notice that $\alpha(y,x)=1-\alpha(x,y)$.
Moreover, $\alpha(x,y)$ is nothing else but the {\it transition probability from $v_{xy}$ to $x$} of the reversible Markov chain associated with $\Gamma^S$.  In additon, if $x\not\sim y$ we define  $\alpha(x,y)=\alpha(y,x)=0.$ Notice that, for any $x,y$  is $\alpha(x,y)=\alpha(y,x)$ iff $c(x,v_{xy})=c(y,v_{xy})=2c(x,y)$.

We  also define, for each $h\in\mathcal{C}(V^S)$ and $u\in\mathcal{C}(V),$  the {\it contraction} of $h$ to $V, \underline{h}\in\mathcal{C}(V),$  as
\begin{equation}
\label{contraction}
\underline{h}(x)=h(x)+\sum\limits_{y\sim x}\alpha(x,y)h(v_{xy}),\hspace{0.5cm} x\in V,\end{equation}
and, the {\it extension of $u$ to $V^S$ with respect to $h$}, $u^h\in\mathcal{C}(V^S),$   as
\begin{equation}
\begin{array}{rll}
\label{extension}
u^h(x)=&\hspace{-0.25cm}u(x)& \hbox{ for all } x\in V;\\[1ex]
u^h(v_{xy})=&\hspace{-0.25cm}\dfrac{h(v_{xy})}{k(v_{xy})} + \alpha(x,y)u(x) + \alpha(y,x)u(y),&  \mbox{ for all } v_{xy}\in V'.
\end{array}\end{equation}
Notice that the extension of $u$ to $V^S$ with respect to $h$, has been defined in order to satisfy  $\L^S(u^h)(v_{xy})=h(v_{xy})$.

\begin{theorem}
Given $h\in \mathcal{C}(V^S)$ such that $\langle h,{\sf 1}_{V^S}\rangle=0,$ then $\langle \underline{h},{\sf 1}_V\rangle=0.$ Moreover, $\overline{u}\in \mathcal{C}(V^S)$ is a solution of the Poisson equation $\mathcal{L}^S(\overline{u}) = h$ in $V^S$ iff $u=\overline{u}_{|V}$ is a solution of the Poisson equation $\mathcal{L}(u) = \underline{h}$ in $V.$ In this case, the identity $\overline{u}=u^{h}$ holds.
\end{theorem}

\demo Firstly we note that $\left\langle \underline{h}, {\sf 1}_V \right\rangle=\left\langle h, {\sf 1}_{V^S} \right\rangle$ as
\vspace{-.85cm}
\begin{center}
$$\sum_{x\in V} \underline{h}(x)=\sum_{x\in V} h(x)+\sum_{x\in V}\sum_{y\sim x} \alpha(x,y)h(v_{xy})=\sum\limits_{x\in V} h(x)+\sum\limits_{v_{xy}\in V'} h(v_{xy}).$$
\end{center}
\vspace{-0.25cm}
So the first statement holds.

Given $h\in \mathcal{C}(V^S)$ such that $\left\langle h, {\sf 1}_{V^S}\right\rangle=0$ and $\overline{u}$ a solution of the Poisson equation $\mathcal{L}^S(\overline{u})=h$ in $V^S,$ then
\vspace{-0.65 cm}
\begin{center}
$$\begin{array}{rlll}
h(v_{xy})=&\hspace{-0.25cm} c(x,v_{xy})\left(\overline{u}(v_{xy})-\overline{u}(x)\right)+c(y,v_{xy})\left(\overline{u}(v_{xy})-\overline{u}(y)\right),&\hbox{for any }v_{xy}\in V';\\[0.25cm]
h(x)=&\hspace{-0.25cm} \sum\limits_{y\sim x} c(x,v_{xy})\left(\overline{u}(x)-\overline{u}(v_{xy})\right),&\hbox{for any }x\in V.
\end{array}$$
\end{center}
\vspace{-0.25cm}
The first identity implies
$\overline{u}(v_{xy})=u^h(v_{xy}),$ assuming $u=\overline{u}_{|V}.$  Then, substituting the expression of $\overline{u}(v_{xy})$ in the second one, we obtain that
\vspace{-1cm}
\begin{center}
$$
\begin{array}{rl}
\mathcal{L}^S(\overline{u})(x)=&\hspace{-0.25cm}\displaystyle\sum\limits_{y\sim x} c(x,v_{xy})\left(\overline{u}(x)-\dfrac{h(v_{xy})}{k(v_{xy})}-\alpha(x,y)\overline{u}(x)-\alpha(y,x)\overline{u}(y)\right)\\[1ex]
=& \hspace{-0.25cm}\displaystyle\sum\limits_{y\sim x} c(x,v_{xy})\alpha(y,x)\left(\overline{u}(x)-\overline{u}(y)\right)-\sum\limits_{y\sim x} \dfrac{c(x,v_{xy})}{k(v_{xy})}h(v_{xy})\\[1ex]
=&\hspace{-0.25cm}\displaystyle\sum\limits_{y\sim x} c(x,y)\left(u(x)-u(y)\right)-\sum\limits_{y\sim x} \alpha(x,y) h(v_{xy})\\[1ex]
=&\hspace{-0.25cm}\mathcal{L}(u)(x)-\underline{h}(x)+h(x),
\end{array}
$$
\end{center}

for every $x\in V.$

Therefore, $\mathcal{L}^S(\overline{u})=h$ in $V^S\,$ iff $\mathcal{L}(u)= \underline{h}$ in $V.$\qed

Next result  shows how to obtain the unique solution of a Poisson problem on the subdivision network $\Gamma^S$
orthogonal to ${\sf 1}_{V^S}.$

\begin{corollary}\label{ortho}
Given $h\in \mathcal{C}(V^S)$, such that $\langle h,{\sf 1}_{V^S}\rangle=0$,  let $\underline{h}\in\mathcal{C}(V)$ be its contraction to $V,$ $u\in \mathcal{C}(V)$ be the unique solution of $\mathcal{L}(u)= \underline{h}$ that satisfies $\langle u,{\sf 1}_V\rangle=0$ and the constant
$$\lambda=-\dfrac{1}{(n+m)}\sum\limits_{x\sim y}\frac{h(v_{xy})}{k(v_{xy})} -\dfrac{1}{(n+m)}\sum\limits_{x\sim y}[ \alpha(x,y)u(x) + \alpha(y,x)u(y)].$$
Then,  $u^\perp=u^h + \lambda$ is the unique solution of $\mathcal{L}^S(u^\perp)= h$ that satisfies $\langle u^\perp,{\sf 1}_{V^S}\rangle=0.$
\end{corollary}

\proof
As two solutions differ on a constant, we have that $u^\perp=u^h+\gamma {\sf 1}_{V^S}$, $\gamma\in\mathbb{R}$. Then,
$$\begin{array}{rl}
0=&\hspace{-.25cm}\displaystyle\langle u^\perp,{\sf 1}_{V^S}\rangle=\langle u^h,{\sf 1}_{V^S}\rangle+(n+m)\gamma=\sum\limits_{x\in V}u(x)+
\sum\limits_{x\sim y}u^h(v_{xy})+(n+m)\gamma\\[2ex]
=&\hspace{-.25cm}\displaystyle\sum\limits_{x\sim y}\frac{h(v_{xy})}{k(v_{xy})} +\sum\limits_{x\sim y}( \alpha(x,y)u(x) + \alpha(y,x)u(y))+(n+m)\gamma,
\end{array}$$
because $\langle u,{\sf 1}_V\rangle=0,$ and the result follows taking $\gamma=\lambda$. \qed

\section{The Green  kernel of a subdivision network}

Taking into account the relation between Poisson problems on $\Gamma^S$ and $\Gamma$, we  obtain the expression of the Green   kernel of a subdivision network, $G^S,$ in terms of   Green's kernel of the base network. From now on we consider the function on $\C(V)$, $\pi^S(x)=\sum\limits_{y\sim x}\alpha(x,y)$ and  the constant $$\beta=\displaystyle\dfrac{1}{(n+m)^2}\sum\limits_{x,y\in V}G(x,y)\pi^S(x)\pi^S(y)
+\dfrac{1}{(n+m)^2}\sum\limits_{x\sim y}\dfrac{1}{k(v_{xy})}.$$

\begin{propo} \label{Gsubdivision}
Let $\Gamma^S$ be the subdivison network of $\Gamma$, then for any $x,z\in V$ and $v_{xy},v_{zt}\in V'$, the Green kernel of $\Gamma^S$ is given by
$$\begin{array}{rl}
G^S(x,z)=&\hspace{-.25cm}\displaystyle G(x,z)-\dfrac{1}{n+m}\sum\limits_{\ell\in V}\Big[G(x,\ell)+G(z,\ell)\Big]\pi^S(\ell)+ \beta,\\[3ex]
 G^S(v_{xy},z)=&\displaystyle\hspace{-.25cm}\alpha(x,y) G(x,z) + \alpha(y,x) G(y,z)\\[3ex]
 -&\hspace{-0.25cm}\displaystyle  \dfrac{1}{n+m}\sum\limits_{\ell\in V}\Big[\alpha(x,y)G(x,\ell)+
 \alpha(y,x)G(y,\ell)+ G(z,\ell)\Big] \pi^S(\ell)-\dfrac{1}{(n+m) k(v_{xy})}+ \beta,\\[3ex]
G^S(v_{xy},v_{zt})=&\hspace{-0.25cm}\displaystyle \alpha(z,t)\Big(\alpha(x,y)G(x,z)+\alpha(y,x)G(y,z)\Big)+\alpha(t,z)\Big(\alpha(x,y)G(x,t)+\alpha(y,x)G(y,t)\Big)\\[3ex]
 -&\hspace{-0.25cm}\displaystyle\dfrac{1}{n+m}\sum\limits_{\ell\in V}\Big[\alpha(x,y)G(x,\ell)+\alpha(y,x)G(y,\ell)+\alpha(z,t)G(z,\ell)+\alpha(t,z)G(t,\ell)\Big]\pi^S(\ell)\\[3ex]
+&\hspace{-0.25cm}\displaystyle\dfrac{\varepsilon_{v_{zt}}(v_{xy})}{k(v_{xy})} -\dfrac{1}{(n+m) k(v_{xy})} -\displaystyle \dfrac{1}{(n+m)k(v_{zt})} +\beta.
\end{array}$$
\end{propo}
\proof
Suppose $z\in V,$ and let $h_z=\varepsilon_z-\dfrac{1}{n+m}.$ Then, for every $x\in V$
$$\underline h_z(x)=\displaystyle\varepsilon_z(x)-\dfrac{1}{n+m}-\dfrac{1}{n+m}\sum\limits_{y\sim x}\alpha(x,y)=\displaystyle\varepsilon_z(x)-\dfrac{1}{n+m}(1+\pi^S(x)).$$
 Hence, from Equation \eqref{Greencharac},  the Poisson problem to solve is $\L(u_z)=\underline h_z,$ and, using the Green kernel for $\Gamma,$ we obtain
$$u_z(x)=G(\displaystyle\varepsilon_z)(x)-\dfrac{1}{n+m}\sum\limits_{\ell\in V}G(x,\ell)\pi^S(\ell)=G(x,z)-\dfrac{1}{n+m}\sum\limits_{\ell\in V}G(x,\ell)\pi^S(\ell).$$
Then, from Corollary \ref{ortho}
$$\begin{array}{rl}
G_z^S(x)=&\hspace{-0.25cm}\displaystyle u_z^{h_z}(x)-\dfrac{1}{(n+m)}\sum\limits_{r\sim s}\frac{h_z(v_{rs})}{k(v_{rs})} -\dfrac{1}{(n+m)}\sum\limits_{r\sim s}[ \alpha(r,s) u_z(r) + \alpha(s,r) u_z(s)]\\[3ex]
=&\hspace{-0.25cm}\displaystyle G(x,z)-\dfrac{1}{n+m}\sum\limits_{\ell\in V}G(x,\ell)\pi^S(\ell)+\dfrac{1}{(n+m)^2}\sum\limits_{r\sim s}\dfrac{1}{k(v_{rs})}\\[3ex]
 -&\hspace{-0.25cm}\displaystyle\dfrac{1}{(n+m)}\sum\limits_{r\sim s}\alpha(r,s)\left[G(r,z)-\dfrac{1}{n+m}\sum\limits_{\ell\in V}G(r,\ell)\pi^S(\ell)\right] \\[3ex]
 -&\hspace{-0.25cm}\displaystyle\dfrac{1}{(n+m)}\sum\limits_{r\sim s}\alpha(s,r)\left[G(s,z)-\dfrac{1}{n+m}\sum\limits_{\ell\in V}G(s,\ell)\pi^S(\ell)\right]\\[3ex]
 =&\hspace{-0.25cm}\displaystyle G(x,z)-\dfrac{1}{n+m}\sum\limits_{\ell\in V}\Big[G(x,\ell)+G(z,\ell)\Big]\pi^S(\ell)+\displaystyle\dfrac{1}{(n+m)^2}\sum\limits_{r, s}G(s,r)\pi^S(r)\pi^S(s)\\[3ex]
+&\dfrac{1}{(n+m)^2}\sum\limits_{r\sim s}\dfrac{1}{k(v_{rs})}.
 \end{array}$$
Now, if
$z\in V$ for every $v_{xy}\in V'$
 $$\begin{array}{rl}
G_z^S(v_{xy})=&\hspace{-0.25cm}\displaystyle \dfrac{h_z(v_{xy})}{k(v_{xy})}
 +\alpha(x,y)u_z(x) + \alpha(y,x)u_z(y)\\[3ex]
 -&\hspace{-0.25cm}\displaystyle\dfrac{1}{(n+m)}\sum\limits_{r\sim s}\frac{h_z(v_{rs})}{k(v_{rs})} -\dfrac{1}{(n+m)}\sum\limits_{r\sim s}[ \alpha(r,s) u_z(r) + \alpha(s,r) u_z(s)]\\[3ex]
 =&\hspace{-0.25cm}\displaystyle -\dfrac{1}{(n+m) k(v_{xy})}
 +\alpha(x,y) G(x,z) + \alpha(y,x) G(y,z)\\[3ex]
 -&\hspace{-0.25cm}\displaystyle  \dfrac{1}{n+m}\sum\limits_{\ell\in V}\big[\alpha(x,y)G(x,\ell)+
 \alpha(y,x)G(y,\ell)\big] \pi^S(\ell)\\[3ex]
 +&\hspace{-0.25cm}\displaystyle\dfrac{1}{(n+m)^2}\sum\limits_{r\sim s}\frac{1}{k(v_{rs})} -\displaystyle\dfrac{1}{(n+m)}\sum\limits_{r\sim s}\alpha(r,s)\left[G(r,z)-\dfrac{1}{n+m}\sum\limits_{\ell\in V}G(r,\ell)\pi^S(\ell)\right] \\[3ex]
 -&\hspace{-0.25cm}\displaystyle\dfrac{1}{(n+m)}\sum\limits_{r\sim s}\alpha(s,r)\left[G(s,z)-\dfrac{1}{n+m}\sum\limits_{\ell\in V}G(s,\ell)\pi^S(\ell)\right]\\[3ex]

=&\hspace{-0.25cm}\displaystyle -\dfrac{1}{(n+m) k(v_{xy})}
 +\alpha(x,y) G(x,z) + \alpha(y,x) G(y,z)\\[3ex]
 -&\hspace{-0.25cm}\displaystyle  \dfrac{1}{n+m}\sum\limits_{\ell\in V}\big[\alpha(x,y)G(x,\ell)+
 \alpha(y,x)G(y,\ell)+ G(z,\ell)\big] \pi^S(\ell)\\[3ex]
+&\hspace{-0.25cm} \displaystyle\dfrac{1}{(n+m)^2}\sum\limits_{r\sim s}\frac{1}{k(v_{rs})}+\dfrac{1}{(n+m)^2}\sum\limits_{r, s}G(s,r)\pi^S(r)\pi^S(s).
 \end{array}$$

Suppose now $v_{zt}\in V,$ and let $h_{v_{zt}}=\varepsilon_{v_{zt}}-\dfrac{1}{n+m}.$ Then, for every $x\in V$
$$\begin{array}{rl}
\underline h_{v_{zt}}(x)=&\hspace{-.25cm}
\displaystyle\varepsilon_{v_{zt}}(x)-\dfrac{1}{n+m}+\sum\limits_{y\in V}\alpha(x,y)\left(\varepsilon_{v_{zt}}(v_{xy})-\dfrac{1}{n+m}\right)\\[3ex]
=&\hspace{-.25cm}
\displaystyle-\dfrac{1}{n+m}(1+\pi^S(x))+\alpha(z,t)\varepsilon_z(x)+\alpha(t,z)\varepsilon_t(x).\end{array}$$
 Hence, the Poisson problem to solve is $\L(u_{v_{zt}})=\underline h_{v_{zt}},$ and, using  Green's kernel for $\Gamma,$ we obtain
$$u_{v_{zt}}(x)=-\dfrac{1}{n+m}\sum\limits_{\ell\in V}G(x,\ell)\pi^S(\ell)+\alpha(z,t)G(x,z)+\alpha(t,z)G(x,t).$$
Then, from Corollary \ref{ortho}
 $$\begin{array}{rl}
G_{v_{zt}}^S(v_{xy})=&\hspace{-0.25cm}\displaystyle \dfrac{h_{v_{zt}}(v_{xy})}{k(v_{xy})}
 +\alpha(x,y)u_{v_{zt}}(x) + \alpha(y,x)u_{v_{zt}}(y)\\[3ex]
 -&\hspace{-0.25cm}\displaystyle\dfrac{1}{(n+m)}\sum\limits_{r\sim s}\frac{h_{v_{zt}}(v_{rs})}{k(v_{rs})} -\dfrac{1}{(n+m)}\sum\limits_{r\sim s}[ \alpha(r,s) u_{v_{zt}}(r) + \alpha(s,r) u_{v_{zt}}(s)]\\[3ex]
 =&\hspace{-0.25cm}\displaystyle \dfrac{\varepsilon_{v_{zt}}(v_{xy})}{k(v_{xy})} -\dfrac{1}{(n+m) k(v_{xy})} \\[3ex]
 -&\hspace{-.25cm}\displaystyle \dfrac{1}{n+m}\sum\limits_{\ell\in V}\Big(\alpha(x,y)G(x,\ell)+\alpha(y,x)G(y,\ell)\Big)\pi^S(\ell)\\[3ex]
 +&\hspace{-.25cm}\displaystyle \alpha(z,t)\Big(\alpha(x,y)G(x,z)+\alpha(y,x)G(y,z)\Big)+\alpha(t,z)\Big(\alpha(x,y)G(x,t)+\alpha(y,x)G(y,t)\Big)\\[3ex]
 -&\hspace{-.25cm}\displaystyle \dfrac{1}{(n+m)k(v_{zt})} +\dfrac{1}{(n+m)^2}\sum\limits_{r\sim s}\frac{1}{k(v_{rs})}\\[3ex]
 -&\hspace{-0.25cm}\displaystyle\dfrac{1}{n+m}\sum\limits_{\ell\in V}\Big(\alpha(z,t)G(z,\ell)+\alpha(t,z)G(t,\ell)\Big)\pi^S(\ell)\\[3ex]
+&\displaystyle \dfrac{1}{(n+m)^2}\sum\limits_{r,s\in V}G(r,s)\pi^S(r)\pi^S(s)\\[3ex]
=&\hspace{-0.25cm}\displaystyle \alpha(z,t)\Big(\alpha(x,y)G(x,z)+\alpha(y,x)G(y,z)\Big)+\alpha(t,z)\Big(\alpha(x,y)G(x,t)+\alpha(y,x)G(y,t)\Big)\\[3ex]
 -&\hspace{-0.25cm}\displaystyle\dfrac{1}{n+m}\sum\limits_{\ell\in V}\Big(\alpha(x,y)G(x,\ell)+\alpha(y,x)G(y,\ell)+\alpha(z,t)G(z,\ell)+\alpha(t,z)G(t,\ell)\Big)\pi^S(\ell)\\[3ex]
 +&\hspace{-0.25cm}\displaystyle\dfrac{\varepsilon_{v_{zt}}(v_{xy})}{k(v_{xy})} -\dfrac{1}{(n+m) k(v_{xy})} -\displaystyle \dfrac{1}{(n+m)k(v_{zt})} +\dfrac{1}{(n+m)^2}\sum\limits_{r\sim s}\frac{1}{k(v_{rs})}\\[3ex]
+&\displaystyle \dfrac{1}{(n+m)^2}\sum\limits_{r,s\in V}G(r,s)\pi^S(r)\pi^S(s).\qed
 \end{array}$$

In particular, if $\Gamma$ is a $k$--regular graph and we consider the standard subdivision graph; that is $c(x,v_{xy})=c(y,v_{xy})=2,$ we get the following result.

\begin{corollary} 
Let $\Gamma^S$ be the standard subdivision graph of  a $k$--regular graph, $\Gamma$; then for any $x,z\in V$ and $v_{xy},v_{zt}\in V'$, the Green kernel of $\Gamma^S$ is given by
$$\begin{array}{rl}
G^S(x,z)=&\hspace{-.25cm}\displaystyle G(x,z)+\dfrac{k}{2n(2+k)^2},\\[3ex]
 G^S(v_{xy},z)=&\displaystyle\hspace{-.25cm}\dfrac{1}{2}\Big( G(x,z) + G(y,z)\Big)-\dfrac{1}{n(2+k)^2},\\[3ex]
G^S(v_{xy},v_{zt})=&\hspace{-0.25cm}\displaystyle \dfrac{1}{4}\Big(G(x,z)+G(y,z)+G(x,t)+G(y,t)+\varepsilon_{v_{zt}}(v_{xy})\Big)-\dfrac{(4+k)}{2n(2+k)^2}.
\end{array}$$
\end{corollary}

\section{Effective Resistances and Kirchhoff Index on subdivision networks}

We are now concerned with the relation between effective resistances in a base network $\Gamma$  and the effective resistances, $R^S$,  in a subdivision network $\Gamma^S.$

\begin{theorem}{\label {teo2}}
Let $\Gamma=(V,E,c)$ be a  network and  $\Gamma^S=(V^S,E^S,c)$ its  subdivision network, then
$$\begin{array}{rl}
R^S(x,y)=&\hspace{-.25cm} R(x,y),\\[2ex]
R^S(x,v_{zt})=&\hspace{-.25cm}\displaystyle \dfrac{1}{k(v_{zt})}+\alpha(z,t)  R(x,z)+\alpha(t,z) R(x,t) -  \alpha(z,t)\alpha(t,z) R(z,t),\\[2ex]
R^S(v_{xy},v_{zt})=&\hspace{-.25cm}\dfrac{1}{k(v_{xy})}+\dfrac{1}{k(v_{zt})}\\[2ex]
 -&\hspace{-.25cm} \alpha(x,y)\alpha(y,x)  R(x,y)- \alpha(z,t)\alpha(t,z)  R(z,t)\\[1ex]
 + &\hspace{-.25cm} \alpha(x,y)\alpha(z,t) R(x,z)+ \alpha(x,y)\alpha(t,z)  R(x,t)\\[1ex]
 +&\hspace{-.25cm}\alpha(z,t)\alpha(y,x) R(y,z)+  \alpha(y,x)\alpha(t,z)R(y,t), \,\,\mbox{ for any }\,\, v_{xy}\not=v_{zt}.
\end{array}$$
\end{theorem}
\proof
The proof is a
 direct  consequence of   Proposition \ref{Gsubdivision} and Identity \eqref{resistancegreen}. Let us do the non--trivial case 2. The case 3, can be proved similarly.

$$\begin{array}{rl}
R^S(x,v_{zt})=&\hspace{-.25cm}G^S(x,x)+G^S(v_{zt},v_{zt})-2G^S(x,v_{zt})\\[2ex]
=&\displaystyle\hspace{-.25cm}G(x,x)-\dfrac{2}{n+m}\sum\limits_{\ell\in V}G(x,\ell)\pi^S(\ell)\\[3ex]
+&\hspace{-0.25cm}\displaystyle \alpha(z,t)\Big(\alpha(z,t)G(z,z)+\alpha(t,z)G(t,z)\Big)+\alpha(t,z)\Big(\alpha(z,t)G(z,t)+\alpha(t,z)G(t,t)\Big)\\[2ex]
 -&\hspace{-0.25cm}\displaystyle\dfrac{2}{n+m}\sum\limits_{\ell\in V}\Big[\alpha(z,t)G(z,\ell)+\alpha(t,z)G(t,\ell)\Big]\pi^S(\ell)+\dfrac{\varepsilon_{v_{zt}}(v_{zt})}{k(v_{zt})} -\dfrac{2}{(n+m) k(v_{zt})}  \\[3ex]
-&\hspace{-0.25cm}\displaystyle2\alpha(z,t) G(z,x) -2 \alpha(t,z) G(t,x)\\[3ex]
 +&\hspace{-0.25cm}\displaystyle  \dfrac{2}{n+m}\sum\limits_{\ell\in V}\Big[\alpha(z,t)G(z,\ell)+
 \alpha(t,z)G(t,\ell)+ G(x,\ell)\Big] \pi^S(\ell)+\dfrac{2}{(n+m) k(v_{zt})}\\[2ex]
=&\displaystyle\hspace{-.25cm}\dfrac{1}{k(v_{zt})}+G(x,x)-2\alpha(z,t) G(z,x) -2 \alpha(t,z) G(t,x)\\[3ex]
+&\hspace{-0.25cm}\displaystyle \alpha(z,t)\Big(\alpha(z,t)G(z,z)+\alpha(t,z)G(t,z)\Big)+\alpha(t,z)\Big(\alpha(z,t)G(z,t)+\alpha(t,z)G(t,t)\Big)\\[2ex]
=&\hspace{-0.25cm}\displaystyle\dfrac{1}{k(v_{zt})}+\alpha(z,t)\Big[G(x,x)+G(z,z)-2G(x,z)\Big]  +\alpha(t,z) \Big[G(x,x)+G(t,t)-2G(x,t)\Big]\\[3ex]
 -&\hspace{-0.25cm}\displaystyle\alpha(t,z)  \alpha(z,t) \Big[G(z,z)+G(t,t)-2G(z,t)\Big],\end{array}$$
and hence, the result follows.\qed

Observe that the effective resistance between vertices of the original network remains unchanged, as expected. In particular for the standard subdivision graph we get the following result, which coincides with the obtained in \cite{ChZh07,SuWaZhBu15,Ya14},  up to   the factor 2 due to our (electrically compatible)--choice of the conductances. 

\begin{corollary}{\label {coro2}}
Let $\Gamma=(V,E,c)$ be a  network and  $\Gamma^S=(V^S,E^S,c)$ its standard subdivision network, then
$$\begin{array}{rl}
R^S(x,y)=&\hspace{-.25cm} R(x,y)\\[2ex]
R^S(x,v_{zt})=&\hspace{-.25cm}\displaystyle \dfrac{1+2 R(x,z)+2 R(x,t) -   R(z,t)}{4}\\[2ex]
R^S(v_{xy},v_{zt})=&\hspace{-.25cm}\dfrac{2-R(x,y)-R(z,t)+R(x,z)+  R(x,t)+ R(y,z)+ R(y,t)}{4}, \,\,\mbox{ for any }\,\, v_{xy}\not=v_{zt}.
\end{array}$$
\end{corollary}

Next we obtain an expression for the Kirchhoff index of the subdivision network, ${\sf k}^S$,  in terms of the Kirchhoff index, ${\sf k}$,  of the base network and other parameters.

\begin{theorem}{\label {kirchhoff}}
Let $\Gamma=(V,E,c)$ be a  network and  $\Gamma^S=(V^S,E^S,c)$ its  subdivision network, then
$$\begin{array}{rl}
{\sf k}^S
=&\hspace{-.25cm}\dfrac{n+m}{n}\displaystyle {\sf k}+(n+m)\displaystyle \sum\limits_{x\in V}G(x,x)\pi^S(x)
-\displaystyle \sum\limits_{x,y\in V}G(x,y)\pi^S(x)\pi^S(y)\\[3ex]
-&\hspace{-.25cm}(n+m)\displaystyle\sum\limits_{x\sim y}\alpha(x,y)\alpha(y,x)R(x,y)
+\displaystyle(n+m-1)\sum\limits_{x\sim y}\dfrac{1}{k(v_{xy})}.\end{array}$$
\end{theorem}
\proof
$$\begin{array}{rl}
{\sf k}^S=&\hspace{-.25cm}(n+m)\displaystyle\sum\limits_{x\in V}G^S(x,x)+(n+m)\sum\limits_{v_{xy}\in V'}G^S(v_{xy},v_{xy})\\[4ex]
=&\hspace{-.25cm}\displaystyle \dfrac{(n+m)}{n}\,{\sf k}-2\sum\limits_{x\in V}\sum\limits_{\ell\in V}G(x,\ell)\pi^S(\ell)\\[3ex]
+&\hspace{-0.25cm}(n+m)\displaystyle \sum\limits_{v_{xy}\in V'}\Big(\alpha(x,y)^2G(x,x)+2\alpha(x,y)\alpha(y,x)G(y,x)+\alpha(y,x)^2G(y,y)\Big)\\[3ex]
 -&\hspace{-0.25cm}\displaystyle 2\sum\limits_{v_{xy}\in V'}\sum\limits_{\ell\in V}\Big[\alpha(x,y)G(x,\ell)+\alpha(y,x)G(y,\ell)\Big]\pi^S(\ell)\\[3ex]
+&\displaystyle\sum\limits_{x,y\in V}G(x,y)\pi^S(x)\pi^S(y) +(n+m-1) \sum\limits_{x\sim y}\dfrac{1}{k(v_{xy})} \\[4ex]
=&\hspace{-.25cm}\displaystyle \dfrac{n+m}{n}\,{\sf k}+(n+m)\sum\limits_{x,y\in V}\Big(\alpha(x,y)^2G(x,x)+\alpha(x,y)\alpha(y,x)G(y,x)\Big)\\[3ex]
-&\displaystyle \sum\limits_{x,y\in V}G(x,y)\pi^S(x)\pi^S(y) +\displaystyle(n+m-1) \sum\limits_{x\sim y}\dfrac{1}{k(v_{xy})}  \\[4ex]
=&\hspace{-.25cm}\displaystyle \dfrac{n+m}{n}\,{\sf k}+(n+m)\displaystyle \sum\limits_{x\in V}G(x,x)\pi^S(x)-(n+m)\sum\limits_{x\sim y}\alpha(x,y)\alpha(y,x)R(y,x)\\[3ex]
-&\displaystyle \sum\limits_{x,y\in V}G(x,y)\pi^S(x)\pi^S(y)+\displaystyle(n+m-1)\sum\limits_{x\sim y}\dfrac{1}{k(v_{xy})}.\qed
\end{array}$$

In particular,  the Kirchhoff index of the standard subdivision graph  has  the following expression which,  coincides with   \cite[Th 3.1] {SuWaZhBu15}. In the case of $k$--regular graph the result coincides with  \cite[Th 3.5] {GaLuLi12}.
 
\begin{corollary}
Let $\Gamma^S$ be the standard subdivision network of  a  graph, $\Gamma$; then
$$\begin{array}{rl}
{\sf k}^S
=&\hspace{-.25cm}\dfrac{n+m}{n}\displaystyle {\sf k}+(n+m)\displaystyle \sum\limits_{x\in V}G(x,x)\pi^S(x)
-\displaystyle \sum\limits_{x,y\in V}G(x,y)\pi^S(x)\pi^S(y)+\dfrac{m^2-n^2+n}{4}.\end{array}$$
In particular, if $\Gamma$ is $k$--regular
$$
{\sf k}^S
=\dfrac{(k+2)^2}{4}{\sf k}+\dfrac{(k^2-4)n^2+4n}{16}.
$$
\end{corollary}
%
%
%

\section{Subdivision network of a wheel}
\begin{figure}
\centering\includegraphics[scale=0.3]{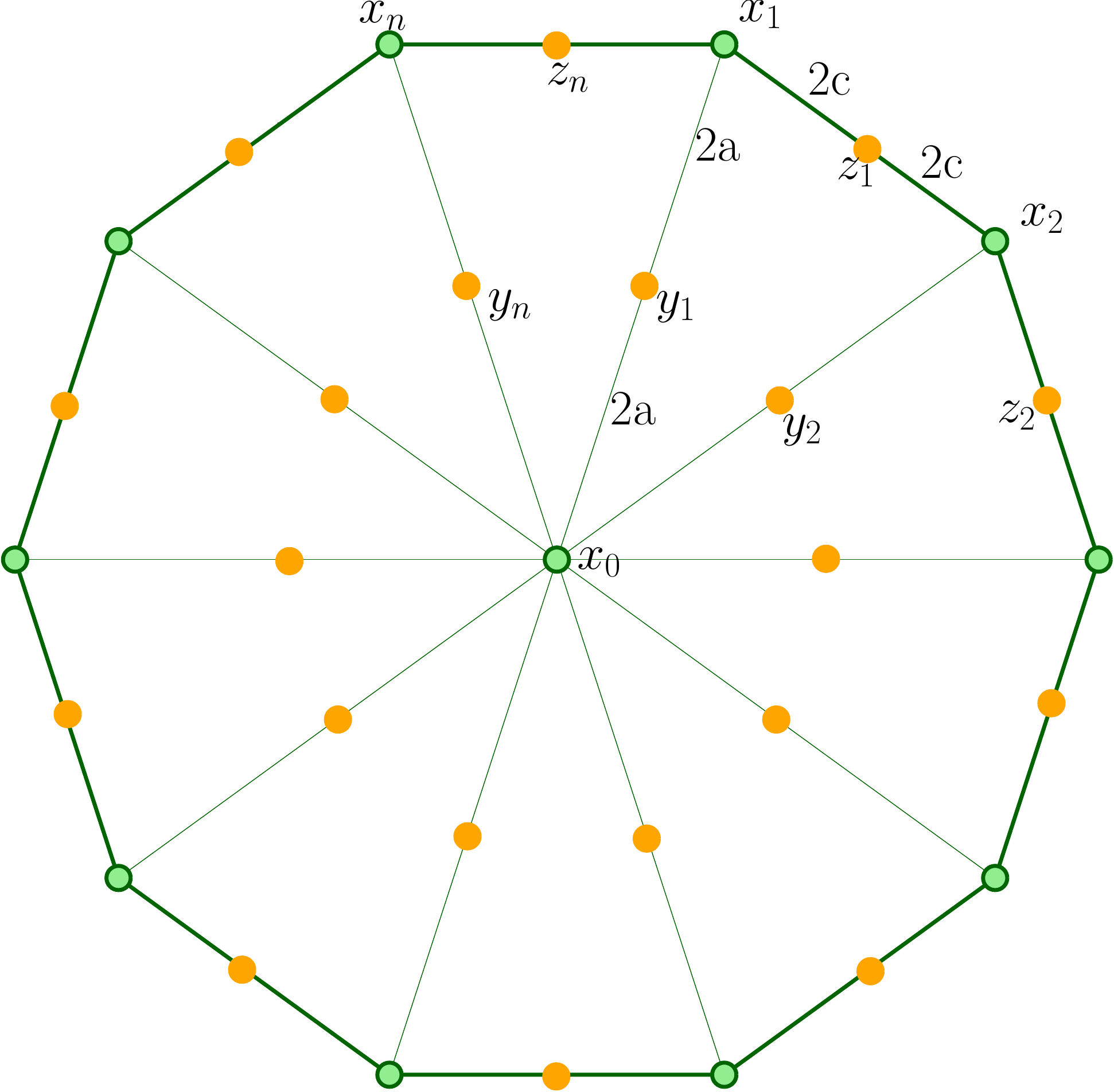}
\caption{\label{wheel} Subdivision network of a wheel of $n$ vertices}
\end{figure}

In order to illustrate the above results, we consider the wheel network with constant conductances and a subdivision of it. Let $W_n$ be the wheel network with vertex set $V=\{x_0, x_1,\ldots,x_n\}$, where $x_0$ has degree $n$, and conductances $c(x_0,x_i)=a>0$ for any $i=1,\ldots, n$,  $c=c(x_i,x_{i+1})$ if $i=1,\ldots,n-1$ and $c=c(x_n,x_1)$, as can be seen in Figure \ref{wheel}. For the sake of simplicity we consider that $x_{n+1}=x_1$.

It is known, see for instance \cite{CaEnMi14-1},
that the  Green  function of $W_n$ is
$$\begin{array}{rl}
G(x_0,x_0)=&\hspace{-.25cm}\dfrac{n}{a(n+1)^2},\\[2ex]
 G(x_0,x_i)=&\hspace{-.25cm}\dfrac{-1}{a(n+1)^2},\quad i=1,\ldots, n,\\[2ex]
G(x_i,x_j)=&\hspace{-.25cm}-\dfrac{n+2}{a(n+1)^2}+\dfrac{U_{n-1-|i-j|}(p)+U_{|i-j|-1}(p)}{2c \big(T_n(p)-1\big)}, \quad i,j=1,\ldots,n,
\end{array}$$
where $p=1+\frac{a}{2c}$ and $U_{\ell}(x), T_{\ell}(x)$ are the Chebyshev polynomials of 1st and 2nd order defined by the recurrence $P_m(x)=2xP_{m-1}(x)-P_{m-2}(x)$ $m\ge 0$ provided that $U_0(x)=1,\ U_1(x)=x$ and $T_{-2}(x)=-1,\ T_{-1}(x)=0,$ respectively.

Let us now define the standard subdivision of the wheel network. The new vertices are $y_i=v_{x_0 x_i}$ and $z_i=v_{x_i x_{i+1}}$ if $i=1,\ldots, n$. The conductances for the new edges  are $2a=c(x_0,y_i)$ and $2c=c(x_i,z_i)$ for $i=1,\ldots, n$. Whereas, the conductance of  the remaining edges  follows taking into account  relation \eqref{conductances}.

Observe that $k( y_i)=4a$ and $k( z_i)=4c$ for $i=1,\ldots,n.$ Moreover, $\alpha(x,y)=\frac{1}{2},$ for every pair of adjacent vertices and $\pi^S(x_0)=\frac{n}{2}$ and $\pi^S(x_i)=\frac{3}{2}, i=1,\ldots,n.$ Then, the expression of the Green kernel for the subdivision network is given next.
\begin{propo} \label{Gsubdivision-wheel}
Let $\Gamma^S$ be the subdivision network of $\Gamma$, and for any $ i,j=1,\ldots,n$ consider $$g_{ij}(p)=\displaystyle \dfrac{U_{n-1-|i-j|}( p )+U_{|i-j|-1}( p )}{2c\big(T_n( p )-1\big)}.$$
Then, the Green kernel for $\Gamma^S$ is given by
$$\begin{array}{rlrl} 
G^S(x_0,x_0)=&\hspace{-.3cm}\displaystyle\dfrac{n(a+26c)}{4ac(3n+1)^2}, & 
G^S(x_0,x_i)=&\hspace{-.3cm}\displaystyle\dfrac{1}{4ac(3n+1)}\left( \dfrac{n(a+26c)}{3n+1}-10c\right),\\[3ex]
G^S(x_0,y_i)=&\hspace{-.3cm}\displaystyle\dfrac{1}{4ac(3n+1)}\Big( \dfrac{n(a+26c)}{3n+1}-6c\Big),
&
G^S(x_0,z_i)=&\hspace{-.3cm}\displaystyle\dfrac{1}{4ac(3n+1)}\Big( \dfrac{n(a+26c)}{3n+1}-(a+10c)\Big),\\[3ex]
\end{array}$$
\vspace{-.25cm}
$$\begin{array}{rl}
G^S(x_i,x_j)=&\hspace{-.25cm}\displaystyle g_{ij}(p)+\dfrac{n(a-34c)-20c}{4ac(3n+1)^2},\\[3ex]
G^S(x_i,y_j)=&\hspace{-0.25cm}\displaystyle \dfrac{1}{2}g_{ij}(p)+\dfrac{n(a-34c)-20c}{4ac(3n+1)^2}+\dfrac{1}{a(3n+1)},\\[3ex]
G^S(x_i,z_j)=&\hspace{-0.25cm}\displaystyle \dfrac{1}{2}\big(g_{ij}(p)+g_{i\,j+1}(p)\big)+\dfrac{n(a-34c)-20c}{4ac(3n+1)^2}-\dfrac{1}{4c(3n+1)},\\[3ex]
G^S(y_i,y_j)=&\hspace{-0.25cm}\displaystyle\dfrac{1}{4}g_{ij}(p)+\dfrac{\varepsilon_{y_i}(y_j)}{4a} +\displaystyle\dfrac{n(a-34c)-20c}{4ac(3n+1)^2}+\dfrac{2}{a(3n+1)},\\[3ex]
G^S(y_i,z_j)=&\hspace{-0.25cm}\displaystyle
\dfrac{1}{4}\big(g_{ij}(p)+g_{i\,j+1}(p)\big)+\dfrac{n(a-34c)-20c}{4ac(3n+1)^2}-\displaystyle \dfrac{a-4c}{4ac(3n+1)}, \\[3ex]
G^S(z_i,z_j) =&\hspace{-0.25cm}\displaystyle\dfrac{p+1}{2}g_{ij}(p)+\displaystyle\dfrac{n(a-34c)-20c}{4ac(3n+1)^2}-\displaystyle\dfrac{1}{2c(3n+1)}+\dfrac{\varepsilon_{z_i}(z_j)}{4c}.
\end{array}$$
\end{propo}
\proof
The expressions given in the proposition follow from the expression for the Green kernel obtained in Proposition \ref{Gsubdivision}. We compute one of the cases in order to illustrate the methodology. 

Firstly, we compute the constant
$$
\beta=\dfrac{1}{(n+m)^2}\sum\limits_{s, r\in V}G(s,r)\pi^S(r)\pi^S(s)
+\dfrac{1}{(n+m)^2}\sum\limits_{r\sim s}\dfrac{1}{k(v_{rs})}=\dfrac{n}{4a(3n+1)^2
}\left[\left(\dfrac{n-3}{n+1}\right)^2+\dfrac{a}{c} +1 \right],
$$
where we have taken into account that
$\sum\limits_{r\in V}G(s,r)=0$ and hence
$$\sum\limits_{ r\in V}G(s,r)\pi^S(r)=\dfrac{n-3}{2}G(s,x_0).$$

Consider $z_i=v_{x_ix_{i+1}}$ and $z_j=v_{x_jx_{j+1}},$ then

$$\begin{array}{rl}
G^S(z_i,z_j)=
&\hspace{-0.25cm}\displaystyle\dfrac{1}{4} \Big(G(x_i,x_j)+G(x_{i+1},x_j)+G(x_i,x_{j+1})+G(
x_{i+1},x_{j+1})\Big)\\[3ex]
 -&\hspace{-0.25cm}\displaystyle\dfrac{1}{2(3n+1)}\sum\limits_{\ell=0}^{n}\Big[G(x_j,x_\ell)+G(x_{j+1},x_\ell)+G(x_j,x_\ell)+G(x_{j+1},x_\ell)\Big]\pi^S(x_\ell)\\[3ex]
+&\hspace{-0.25cm}\displaystyle\dfrac{\varepsilon_{z_j}(z_i)}{k(z_i)} -\dfrac{1}{(3n+1) k(z_i)} -\displaystyle \dfrac{1}{(3n+1)k(z_j)} +\beta\\[3ex]
=&\hspace{-0.25cm}\displaystyle-\dfrac{n+2}{a(n+1)^2} +\dfrac{2U_{|i-j|-1}( p )+U_{|i+1-j|-1}( p )+U_{|i-j-1|-1}( p )}{8c\big(T_n( p )-1\big)}\\[3ex]
+&\hspace{-0.25cm}\dfrac{2U_{n-1-|i-j|}( p )+U_{n-1-|i+1-j|}( p )+U_{n-1-|i-j-1|}( p )}{8c\big(T_n( p )-1\big)}\\[3ex]
 +&\hspace{-0.25cm}\displaystyle\dfrac{n-3}{(3n+1)(n+1)^2}+\displaystyle\dfrac{\varepsilon_{z_j}(z_i)}{4c} -\dfrac{2}{(3n+1) 4c} +\beta\\[3ex]
=&\hspace{-0.25cm}\displaystyle\dfrac{(a+4c)\Big(U_{n-1-|i-j|}( p )+U_{|i-j|-1}( p )\Big)}{8c^2\big(T_n( p )-1\big)}-\dfrac{n(5a+34c)+2a+20c}{4ac(3n+1)^2} +\dfrac{\varepsilon_{z_i}(z_j)}{4c}.\qed
\end{array}$$

%

%

To end up the section we compute the Kirchhoff index of the standard subdivision graph associated with the wheel  $W_n$.
\begin{corollary} The Kirchhoff index of the standard subdivision network os $W_n$ is
$${\sf k}^S=\dfrac{3n^2(a+c)-25cn}{4ac}+\dfrac{n(3n+1)\big(7U_{n-1}(p)+2U_{n-2}(p)+2\big)}{8c\big(T_n( p )-1\big)}.$$
\end{corollary}
\section*{Acknowledgements} This work has been partly supported by the Spanish Research Council (Comisi\'on Interministerial de Ciencia y Tecnolog\'{\i}a,) MTM2014-60450-R.


\begin{thebibliography}{10}

\bibitem{BeCaEnGe08a}
E.~Bendito, A.~Carmona, A.~M. Encinas, and J.~M. Gesto.
\newblock A formula for the {K}irchhoff index.
\newblock {\em Int. J. Quantum Chem.}, 108(6):1200--1206, 2008.

\bibitem{BuYaZhZh14}
C. Bu, B. Yan, X. Zhou, and J. Zhou.
\newblock Resistance distance in subdivision-vertex join and subdivision-edge
  join of graphs.
\newblock {\em Linear Algebra Appl.}, 458:454--462, 2014.

\bibitem{CaEnMi14-1}
A.~Carmona, A.~M. Encinas, and M.~Mitjana.
\newblock Discrete elliptic operators and their {G}reen operators.
\newblock {\em Linear Algebra Appl.}, 442:115--134, 2014.

\bibitem{Che10}
H. Chen.
\newblock Random walks and the effective resistance sum rules.
\newblock {\em Discrete Appl. Math.}, 158(15):1691--1700, 2010.

\bibitem{ChZh07}
H. Chen and F. Zhang.
\newblock Resistance distance and the normalized {L}aplacian spectrum.
\newblock {\em Discrete Appl. Math.}, 155(5):654--661, 2007.

\bibitem{GaLuLi12}
X. Gao, Y. Luo, and W. Liu.
\newblock Kirchhoff index in line, subdivision and total graphs of a regular
  graph.
\newblock {\em Discrete Appl. Math.}, 160(4-5):560--565, 2012.

\bibitem{GhBoSa08}
A.~Ghosh, S.~Boyd, and A.~Saberi.
\newblock Minimizing effective resistance of a graph.
\newblock {\em SIAM Review}, 50(1):37--66, 2008.

\bibitem{KlRa93}
D.J. Klein and M.~Randi\'c.
\newblock Resistance distance.
\newblock {\em J. Math. Chem.}, 12(1):81--95, 1993.


\bibitem{SuWaZhBu15}
L. Sun, W. Wang, J. Zhou, and C. Bu.
\newblock Some results on resistance distances and resistance matrices.
\newblock {\em Linear Multilinear Algebra}, 63(3):523--533, 2015.

\bibitem{XiGu03}
W.~Xiao and I.~Gutman.
\newblock Resistance distance and {L}aplacian spectrum.
\newblock {\em Theor. Chem. Acc.}, 110(4):284--289, 2003.

\bibitem{Ya14}
Y. Yang.
\newblock The {K}irchhoff index of subdivisions of graphs.
\newblock {\em Discrete Appl. Math.}, 171(0):153 -- 157, 2014.

\end{thebibliography}
\end{document}